\documentclass[arxiv]{agn_article}

\usepackage{agn_all}
\usepackage{graphicx}
\usepackage{amsmath,mathtools,amssymb}
\usepackage{agn_macros}
\usepackage{float}
\usepackage{yfonts}
\usepackage{longtable}
\usepackage[hyperindex,breaklinks]{hyperref}
\newcommand{\citet}[2][]{\citeauthor{#2} \cite[#1]{#2}}
\newcommand{\Wiso}{W_{\mathrm{iso}}}
\newcommand{\Wvol}{W_{\mathrm{vol}}}

\pgfplotsset{
every axis/.append style={width=0.6\textwidth, height=6cm},
every axis plot/.append style={no markers, thick},
label style={font=\small},
tick label style={font=\small},
legend pos={outer north east},
cycle list={red, green, cyan, yellow}
}
\numberwithin{equation}{section}

\begin{document}

%language
\selectlanguage{\english}

%macros
\newcommand{\ex}[1]{e^{#1}}
\renewcommand{\.}{\hspace*{0.07em}}

% title page
\title{A commented translation of Hans Richter's early work\\\enquote{The isotropic law of elasticity}\\ 
}
\author{%
	Kai Graban\thanks{% 
		Kai Graban, \quad Fakult\"{a}t f\"{u}r Mathematik, Universit\"{a}t Duisburg-Essen, Thea-Leymann Str. 9, 45127 Essen%
	}\quad\ and\quad%	
	Eva Schweickert\thanks{%
		Eva Schweickert,\quad Lehrstuhl f\"{u}r Nichtlineare Analysis und Modellierung, Fakult\"{a}t f\"{u}r Mathematik, Universit\"{a}t Duisburg-Essen, Thea-Leymann Str. 9, 45127 Essen, Germany; email: eva.schweickert@stud.uni-due.de%
	}\quad\ and\quad%
	Robert J.\ Martin\thanks{%
		Robert J.\ Martin,\quad Lehrstuhl f\"{u}r Nichtlineare Analysis und Modellierung, Fakult\"{a}t f\"{u}r Mathematik, Universit\"{a}t Duisburg-Essen, Thea-Leymann Str. 9, 45127 Essen, Germany; email: robert.martin@uni-due.de%
	}\quad\ and\quad%
	Patrizio Neff\thanks{%
		Patrizio Neff,\quad Head of Lehrstuhl f\"{u}r Nichtlineare Analysis und Modellierung, Fakult\"{a}t f\"{u}r	Mathematik, Universit\"{a}t Duisburg-Essen, Thea-Leymann Str. 9, 45127 Essen, Germany, email: patrizio.neff@uni-due.de%
	}%
}
\date{\today}

\maketitle

\begin{abstract}
	We provide a faithful translation of Hans Richter's important 1948 paper \enquote{Das isotrope Elastizitätsgesetz} from its original German version into English. Our introduction summarizes Richter's achievements.	
\end{abstract}

\\[.65em]

%%%%%%%%%%%%%%%%%%%%%%%%%%%%%%%%%%%%%%%%%%%%%%%%%%%%%
\section*{Introduction}
Shortly after the second world war, in a series of papers \cite{richter1948,richter1949verzerrung,richter1949log,richter1952elastizitatstheorie} from $1948-1952$, Hans Richter ($1912-1978$) laid down his general format of isotropic nonlinear elasticity based on a rather modern approach with direct tensor notation. By translating his work \enquote{Das isotrope Elastizitätsgesetz} \cite{richter1948}, we aim at making his development, which precedes later work in the field by several decades, accessible to the international audience.\\
\noindent\hspace*{0.5cm} Let us briefly summarize  Richter's achievement in this paper. He uses, for the time, rather advanced methods of matrix analysis (including the theory of primary matrix functions \cite{agn_martin2015some}) and employs the left polar decomposition \cite{agn_lankeit2014minimization,agn_neff2014logarithmic,agn_neff2014grioli} of the deformation gradient $F=V\,R$ into a stretch $V\in\Sym^+(3)$ and a rotation $R\in\SO(3)$. For Richter, the \enquote{physical stress tensor} is the Cauchy stress tensor $\sigma\in\Sym(3)$. From the coaxiality between $\sigma$ and $V$ for an isotropic response, he deduces the representation formula for isotropic tensor functions (the \emph{Richter representation}, see \eqref{eq:repformula})
\begin{align}\label{eq:sigma}
	\sigma=g_1(I_1,I_2,I_3)\cdot\id+g_2(I_1,I_2,I_3) \cdot V+g_3(I_1,I_2,I_3)\cdot V^2\,	
\end{align} 
(predating the Rivlin-Ericksen representation theorem \cite{rivlin1955} by $7$ years) where $g_i$, $i=1,2,3$ are scalar valued functions of the invariants $I_\nu$, $\nu=1,2,3$, with
\begin{align*}
	I_1=\tr(V),\quad I_2=\frac{1}{2}\tr(V^2), \quad I_3=\det V.
\end{align*} 
Alongside, Richter introduces the logarithmic stretch tensor $L=\log V$ without citing the previous work of Hencky \cite{Hencky1928,Hencky1929,hencky1929super,hencky1931,neff2016geometry,agn_neff2015exponentiatedI,agn_neff2015exponentiatedII,agn_neff2014exponentiatedIII,agn_neff2014rediscovering}. He then turns to the question of what happens if the relation \eqref{eq:sigma} is derived from a stored energy $W(I_1,I_2,I_3)$, i.e.\ when \eqref{eq:sigma} is consistent with hyperelasticity. He obtains the correct representation (see \eqref{eq:sigmaW} in his text) 
\begin{align}\label{eq:sigmaW1}
	\sigma=\pdd{W}{I_3}\cdot\id+\tel{I_3}\cdot\pdd{W}{I_1} \cdot V+\tel{I_3}\cdot\pdd{W}{I_2}\cdot V^2\,,\qquad W=W(I_1, I_2, I_3).	
\end{align}
In the next section, Richter introduces the \emph{multiplicative split} of the elastic stretch $V$ into volume preserving (isochoric) parts and volume change (see \eqref{4.1})
\begin{align}
	V=\frac{V}{(\det V)^\frac{1}{3}}\cdot(\det V)^\frac{1}{3}\cdot\id
\end{align}
and he observes that the logarithmic stretch tensor additively separates both contributions by using the classical deviator operation (see \eqref{4.2}) such that 
\begin{align}
	\log V=\dev\log V+\frac{1}{3}\tr(\log V)\cdot\id=\log\left(\frac{V}{(\det V)^\frac{1}{3}}\right)+\frac{1}{3}\log\det V\cdot\id,\quad\dev X=X-\frac{1}{3}\tr(X)\cdot\id.
\end{align}
He also observes that the invariants based on the logarithmic stretch tensor satisfy certain algebraic relations, cf.\ \cite{criscione2000invariant}.
In Richter's fifth section, he introduces the \emph{volumetric-isochoric split}
\begin{align}
	W(F)&=\Wiso\left(\dev\log V\right)+\Wvol\left(\tr(\log V)\right)\notag\\
	&=\Wiso\left(\log\left(\frac{V}{(\det V)^\frac{1}{3}}\right)\right)+\Wvol\left(\log\det V\right)\notag
	=\widetilde{W}_\text{iso}\left(\frac{V}{(\det V)^\frac{1}{3}}\right)+\widetilde{W}_\text{vol}(\det V)
\end{align}
of the stored energy (often erroneously attributed to \cite{flory1961thermodynamic}) and he immediately obtains the important result:\\
\noindent\hspace*{0.5cm} \emph{An isotropic energy is additively split into volumetric and isochoric parts if and only if the mean Cauchy stress $\frac{1}{3}\tr\sigma$ is only a function of the relative volume change $\det V$. In that case,} 
\begin{align}
	\frac{1}{3}\tr\sigma=\frac{1}{\det V}\cdot \Wvol^\prime(\log\det V)=\widetilde{W}^\prime_\text{vol}(\det V).
\end{align}
\noindent\hspace*{0.5cm} This result has been rediscovered and re-derived multiple times, e.g.\ in \cite{charrier1988existence,murphy2018modelling,sansour2008physical,favrie2014thermodynamically,federico2010volumetric,ndanou2017piston}. In addition, Richter shows that this property of the volumetric-isochoric split is invariant under a change of the reference temperature. Finally, he poses the question whether a linear relation between $\sigma$ and $V$ in the form (the Hooke's law as he perceives it)
\begin{align}\label{eq:sigmalinear}
	\sigma=2\mu\,(V-\id)+\lambda\tr(V-\id)\cdot\id,
\end{align}
where $\mu>0$ is the shear modulus and $\lambda$ is the second Lamé parameter, can be consistent with hyperelasticity. A short calculus reveals that \eqref{eq:sigmalinear} is hyperelastic if and only if $2\mu=\lambda$, i.e.\ for Poisson ratio $\nu=\frac{1}{3}$ (which is approximately satisfied for many metals, e.g.\ aluminium). For all other values of $\nu$, Hooke's law is incompatible with the hyperelastic approach and Richter proposes to use instead (the quadratic Hencky energy \cite{Hencky1928,agn_neff2014axiomatic})
\begin{align*}
	W(F)=\mu\,\norm{\dev\log V}^2+\frac{2\mu+3\lambda}{6}\,\tr^2(\log V)
\end{align*}
with the induced stress-strain law
\begin{align}
	\sigma\cdot\det F=\tau=2\mu\,\log V+\lambda\,\tr(\log V)\cdot\id\,,
\end{align}
where $\tau$ is the Kirchhoff stress tensor.\\
\noindent\hspace*{0.5cm} We will briefly discuss the constitutive relation \eqref{eq:sigmalinear}. In order to check hyperlasticity of the Cauchy stress-stretch relation in this case,
we use the representation, consistent with \eqref{eq:sigmaW1},
\begin{align}
	\sigma(V)=\frac{2\,\mu}{J}\, D_V W(V)\cdot V, \qquad J=\det V,
\end{align}
and consider the energy $W(F)=2\mu\,\det V\, [\tr(V)-4]$. Then $\sigma(V)=2\mu\,(V-\id)+2\mu\,\tr(V-\id)\cdot\id$.

Since $\tr[\sigma(V)]=8\mu\,\tr(V-\id)$ and $\tr[\sigma(\alpha \,\id)]=24\mu\,(\alpha-1)$, the Cauchy stress tensor given by \eqref{eq:sigmalinear} with $2\mu=\lambda$ is injective (but not bijective, since for $\tr(\sigma)=-K^+<-24\mu$ there does not exist a stretch $V\in\PSym(3)$ such that $\tr(\sigma(V))=-K^+$).
Furthermore, note that
\begin{align}
	[\tr (V)]^2=(\lambda_1+\lambda_2+\lambda_3)^2=\lambda_1^2+\lambda_2^2+\lambda_3^2+2(\lambda_1\lambda_2+\lambda_1\lambda_3+\lambda_2\lambda_3)=\tr(B)+2\tr(\Cof B),
\end{align}
where $\lambda_i$ are the singular values of the deformation gradient $F$. Then
\begin{align}
	2\mu\,\det V\, [\tr(V)-4]&=2\mu\,\sqrt{\det(B)}\,\{\sqrt{\tr(B)+2\tr(\Cof B)}-4\}
	=2\mu\,\sqrt{I_3}\,\{\sqrt{I_1+2\,I_2}-4\}\notag\\
	&=W(I_1,I_2,I_3)\,,\qquad I_1=\tr(B),\quad I_2=\tr(\Cof B),\quad I_3=\det B.
\end{align}
For this energy, the weak empirical inequalities \cite{agn_thiel2019empirical} $\frac{\partial W}{\partial I_1}>0$ and $\frac{\partial W}{\partial I_2}>0$ are satisfied.
The principal Cauchy stresses are given by $\sigma_i=2\mu\cdot\left(\lambda_i-1+(\lambda_1+\lambda_2+\lambda_3-3)\right)$, which shows that the tension-extension (TE) inequalities and the Baker-Ericksen (BE) inequalities \cite{bakerEri54}, given by
\[
	0 < \frac{\partial \sigma_i}{\partial \lambda_i}=2\mu\cdot\left(1+1\right)=4\mu
	\qquad\text{and}\qquad
	0 < (\sigma_i-\sigma_j)(\lambda_i-\lambda_j)=2\mu\,(\lambda_i-\lambda_j)^2
\]
respectively, are satisfied as well.
We also note that $W(V)=2\mu\cdot\det V\cdot [\tr(V)-4]$ is the Shield-transformation \cite{shield1967inverse} of $W^*(F)=2\mu\cdot [\tr(V^{-1})-4]$, where
\begin{align}
	W^*(F)&=2\mu\left(\frac{1}{\lambda_1}+\frac{1}{\lambda_2}+\frac{1}{\lambda_3}-4\right)=g(\lambda_1,\lambda_2,\lambda_3)
\end{align}
has the Valanis-Landel form\footnote{This calculus shows that the Valanis-Landel form is not invariant under the Shield transformation. In addition, the mapping $U\mapsto T_{\rm Biot}=D_UW(U)$ of the stretch $U=\sqrt{F^TF}$ to the Biot stress tensor $T_{\rm Biot}$ is strictly monotone.} \cite{valanis1967} and $g$ is convex in $(\lambda_1,\lambda_2,\lambda_3)$; the TE-inequalities are satisfied as well.

\medskip
\noindent\hspace*{0.5cm}Richter's paper is not only written in German, but his notation strongly relies on German fraktur letters, which makes reading his original work rather challenging.
In our faithful translation of his paper, we have therefore updated the notation to more current conventions; a complete list of notational changes is provided in Appendix \ref{appendix:listOfSymbols}. Richter's original equation numbering has been maintained throughout.

{\footnotesize \printbibliography}\newpage
%%%%%%%%%%%%%%%%%%%%%%%%%%%%%%%%%%%%%%%%%%%%%%%%%%%%%%%%%%%
\begin{center}
	{ \huge The isotropic law of elasticity\\\medskip \large By \emph{Hans Richter} in Haltingen (Lörrach)\\[0.5cm]
	Zeitschrift für Angewandte Mathematik und Mechanik, Vol. $28$, $1948$, page $205-209$}
\end{center} 
\begin{abstract}
	From the demand of the isotropy and of the existence of the thermodynamic potentials a  general form of the three-dimensional law of elasticity is stated. In doing so, the logarithmic matrix of relative elongations is used, which permits the separation of the variation of the volume and that of the shape by simply forming the deviator. The resilience energy is exactly the sum of the energy of the variation of the volume and that of the shape, if the average tension depends only on the variation of the volume. For finite deformations, the law of \emph{Hooke} is permissible only in the case $\nu=\frac{1}{3}$. 
	\par \medskip \noindent
	Aus der Forderung der Isotropie und der Existenz der thermodynamischen Potentiale wird für das räumliche Elastizitätsgesetz eine allgemeine Form angegeben, wobei die logarithmische Dehnungsmatrix verwendet wird, bei der die Trennung in Volum- und Gestaltänderung durch gewöhnliche Deviatorbildung möglich ist. Die elastische Energie ist genau dann die Summe aus Volum- und Gestaltänderungsenergie, wenn die mittlere Spannung nur von der Volumänderung abhängt. Das \emph{Hookesche} Gesetz ist für endliche Verzerrungen nur bei $\nu=\frac{1}{3}$ zulässig.
	\par \medskip \noindent
	En supposant l'isotropie et l'existence des potentiels thermodynamiques, on donne une forme g\'{e}n\'{e}rale de la loi de l'\'{e}lasticit\'{e} en se servant d'une matrix logarithmique d'allongement. Ce proc\'{e}d\'{e} permet une s\'{e}paration des changements de volume et de forme par une simple formation de d\'{e}viateur. Si la tension moyenne ne d\'{e}pend que du changement de volume, l'\'{e}nergie d'\'{e}lasticit\'{e} est la somme des \'{e}nergies de changement du volume et de la forme. La loi de \emph{Hooke} n'est admissible que pour $\nu=\frac{1}{3}$.
\end{abstract}

%%%%%%%%%%%%%%%%%%%%%%%%%%%%%%%%%%%%%%%%%%%%%%%%
\section{Definitions}
In generalization of Hooke's law, a material is called purely elastic if the Cauchy stresses depend in a uniquely reversible way on the stretches. Strictly speaking, however, it is necessary to discuss the heat transfer which occurs in the tensile test; in particular, it is necessary to distinguish between an adiabatic and isothermal law of elasticity. This choice also clarifies what is meant by strains, since strains on the adiabat resp.\ isotherm can be referred e.g.\ to the initial state, for which the stresses disappear completely. The strains can also be referred to a stress-free initial state at an arbitrarily chosen initial temperature $\Theta_0$ instead. Then the stress-free state at another temperature $\Theta$ corresponds, in the case of an isotropic material, to uniform stretches in all directions, i.e.\ the thermal expansion. In this manner the law of thermal expansion is included in the elastic law. Of course, the affected material must be assumed not to change permanently by changes in temperature within the considered temperature range.
\par \noindent
Thus, we assume a stress-free state at a temperature $\Theta_0$. Let the deformation of the material into another state be characterized by the matrix $F$ and the related stresses by the stress tensor $\sigma$.\footnote{$F$ is the Jacobian matrix: $\intd{\widehat{x}}=F\.\intd{x}$. $\sigma$ is the physical stress tensor at the point $\widehat{x}$.} We call the material ideally elastic if $\sigma$ depends uniquely on $F$ and $\Theta$. The material is said to be isotropic if this dependence is invariant under Euclidean rotations.
\par \noindent
When solving the problem of finding the most general form of this dependence, one appropriately operates with matrices, where the following abbreviations are used:
\par \medskip \noindent
$X^T$ is the matrix obtained by reflecting $X$ over its main diagonal. $(X)_{ik}$ is the entry in the $i$-th row and the $k$-th column of $X$. $\det X$ is the determinant of $X$. $\tr X$ is the sum of the elements on the main diagonal of $X$: called the trace of $X$. $\id$ is the identity tensor. If $f(x)=\sum a_n \cdot x^n$, then, assuming convergence, $f(X)=\sum a_n\.X^n$.
\par \noindent
Recall the following simple statements:
\begin{align}
	\tr(X\cdot Y)&=\tr(Y\cdot X)\,.\label{1.1}\\
	\tr(X\cdot\intd{\log Y})&=\tr(X\cdot Y^{-1}\cdot\intd{Y})\label{1.2}
\end{align}
if $X$ commutes with $Y$, but not necessarily with $\intd{Y}$. 
\begin{align}
	\log(\det X)=\tr(\log X)\,,\label{1.3}
\end{align}
if $\log X$ is well defined.
\begin{alignat}{2}
	&\text{For a pure rotation $R$:}\hspace*{3.5em} & R\.R^T&=\id\,.\hspace*{17.5em} \label{1.4}\\
	&\text{For a pure stretch $V$:} & V&=V^T\,. \notag
\end{alignat}
Every $X$ can be represented in the form:
\begin{align}
	X=V\cdot R\,,\label{1.5}
\end{align}
where the multiplication is to be read in its functional notation from right to left.
%%%%%%%%%%%%%%%%%%%%%%%%%%%%%%%%%%%%%%%%%%%%%%%%%%%%%%%%%%%%%%%%%%%%%%
\section{Consequence of isotropy}
According to \eqref{1.5}, $F$ can be interpreted as a rotation $R$ followed by a stretch $V$, where the principal stretch directions of the latter are rotated against those of the coordinate axes. For the case of isotropic materials, the application of $R$ must not have any influence on $\sigma$. Therefore, $\sigma$ is a function of $V$ and $\Theta$.
For given $F$, we can find $V$ by using \eqref{1.4} and \eqref{1.5} by
\begin{align}
	F\.F^T=V\.R\.R^TV^T=V^2\,.\label{2.1}
\end{align}
The most general coaxial relation between $\sigma$ and $V$ which fulfills the invariance under rotations is now, obviously,
\begin{align}
	\sigma=f(V; I_1, I_2, I_3, \Theta)\,,\label{2.2}
\end{align}
where the $I_\nu$ are the invariants\footnote{It is easy to see that here, one of the invariants $I_ \nu$ can be omitted, in contrast to the subsequent formula \eqref{2.7}.} of $V$.
\par \noindent
Instead of $V$, one can also use a uniquely invertible function of $V$. As we will see later on, it is appropriate to use the ``logarithmic stretch''
\begin{align}
	L=\log V\,,\label{2.3}
\end{align}
which is always defined because of the positive eigenvalues of $V$. We denote the invariants of $L$ by
\begin{align}
j=\tr(L)\,,\qquad k=\tr(L^2)\qquad\text{and}\qquad l=\tr(L^3)\,.\label{2.4}
\end{align}
Further, from \eqref{1.3} and \eqref{2.1} we obtain: \quad $j=\dfrac{1}{2}\.\tr(\log (F\.F^T))=\dfrac{1}{2}\.\log(\det(F\.F^T))=\log(\det F)$\,.
\par \medskip \noindent
Instead of \eqref{2.2}, we can now write
\begin{align}
	\sigma=f(L; j, k, l, \Theta)\,.\label{2.5}
\end{align}
Here, $\tr(\sigma)$, $\tr(\sigma\.L)$ and $\tr(\sigma\.L^2)$ are functions of $j$, $k$, $l$ and $\Theta$ due to \eqref{2.5}.
If we now define the invariants $f_1$, $f_2$ and $f_3$ as the solutions to the system of equations
\begin{align*}
	\tr(\sigma)&=f_1 \tr(\id)+f_2 \tr(L)+f_3 \tr(L^2) \\
	\tr(\sigma\.L)&=f_1 \tr(L)+f_2 \tr(L^2)+f_3 \tr(L^3) \\
	\tr(\sigma\.L^2)&=f_1 \tr(L^2)+f_2 \tr(L^3)+f_3 \tr(L^4)
\end{align*}
with, in general, non-vanishing determinant, then we have for
\begin{align*}
	X=f_1 \cdot \id+f_2 \cdot L+f_3 \cdot L^2\,: \qquad \tr(\sigma\.L^\nu)=\tr(X\.L^\nu) \quad \text{with} \quad \nu=0, 1, 2\,.
\end{align*}
Since $\sigma$ is coaxial to $L$, it is completely determined by $\tr(\sigma)$, $\tr(\sigma\.L)$ and $\tr(\sigma\.L^2)$. Therefore, $\sigma\equiv X$ holds; i.e.
\begin{align}\label{eq:repformula}
	\sigma=f_1(j, k, l, \Theta) \cdot \id+f_2(j, k, l, \Theta) \cdot L+f_3(j, k, l, \Theta)\cdot L^2\,.
\end{align}
Hence, we have found the most general isotropic relation. Using $V$ instead of $L$, we would correspondingly obtain:
\begin{align}
	\sigma=g_1(I_\nu, \Theta)\cdot\id+g_2(I_\nu, \Theta) \cdot V+g_3(I_\nu, \Theta)\cdot V^2\,.\label{2.7}
\end{align}
%%%%%%%%%%%%%%%%%%%%%%%%%%%%%%%%%%%%%%%%%%%%%%%%%%%%%%%%%%%%%%%%%%%%%%%%
\section{Consequence of the	potential}
The internal energy of the material per unit volume in the initial state is denoted by
\begin{align}
	E&=E(j,k,l,\Theta)\,;\label{3.1}\\
\intertext{the entropy is denoted by}
	S&=S(j,k,l,\Theta)\,.\label{3.2}
\end{align}
Then the free energy $W$ takes the form
\begin{align}
	W=E-\Theta\cdot S=W(j,k,l,\Theta)\,.\label{3.3}
\end{align}
If $\intd{A}$ is now the differential of the work done by the element of volume, then
\begin{align}
	\intd{A}=-\intd{E}+\Theta\cdot \intd{S}=-\intd{W}-S\cdot \intd{\Theta}\,.\label{3.4}
\end{align}
Thus for isothermal elastic changes, we have
\begin{align}
	\intd{A}&=-(\intd{W})_{\Theta=\mathrm{const.}}\,;\label{3.5}\\
\intertext{whereas for adiabatic changes}
	\intd{A}&=-(\intd{E})_{S=\mathrm{const.}}\,,\label{3.6}
\end{align}
where $\Theta$ has to be eliminated in \eqref{3.1} and \eqref{3.2}, so that $E$ appears as a function of $j$, $k$, $l$ and $S$.
\par \noindent
In order to calculate $\intd{A}$, we transition from  a deformation $F$ to the neighboring deformation $F+\intd{F}$. Since a pure rotation has no influence on $\intd{A}$, we can assume that $F$ is a pure stretch. Let $e_1$, $e_2$ and $e_3$ be the unit vectors in the principal stretch directions of $V$, which can be interpreted as coordinate vectors. Let $\sigma_1$, $\sigma_2$ and $\sigma_3$ be the components of $\sigma$ in these directions. We can use the rectangular parallelepiped spanned by $V\.e_1$, $V\.e_2$ and $V\.e_3$ as the volume element, which is generated by the stretch $V$ applied to the unit cube. Let us now consider the side which starts from $V\.e_1$ and which is spanned by $V\.e_2$ and $V\.e_3$. Besides an infinitesimal tilting and change of the surface, this side undergoes a displacement in the $e_1$-direction with the magnitude $e_1 \cdot ((V+\intd{F})\.e_1-V\.e_1)=e_1\.\intd{F}\.e_1=(\intd{F})_{11}$ in the transition from $V$ to $V+\intd{F}$. The work done on the considered side is therefore
\begin{align*}
	-\sigma_1\cdot(\intd{F})_{11}\cdot (V)_{22}\cdot (V)_{33}=-\det(V)\cdot\frac{(\intd{F})_{11}\cdot\sigma_1}{(V)_{11}}\,.
\end{align*}
Thus the entire work done by the volume element is
\begin{align}
	\intd{A}=-\det(V)\cdot\sum_{v=1}^3\frac{(\intd{F})_{vv}\cdot\sigma_v}{(V)_{vv}}=-\det(V)\cdot\tr(\sigma\.V^{-1}\intd{F})\,.\label{3.7}
\end{align}
The deformation $V+\intd{F}$ now corresponds to a stretch $V+\intd{V}$, where due to \eqref{2.1},
\begin{align*}
	(V+\intd{V})^2=(V+\intd{F})(V+&\intd{F}^T)
\intertext{or}
	V\cdot\intd{V}+\intd{V}\cdot V=V\cdot\intd{F}^T+&\intd{F}\cdot V\,.
\end{align*}
Multiplying the left side of the equation by $\sigma\.V^{-2}$, taking the trace and using \eqref{1.1}, we find
\begin{align*}
	2\.\tr(\sigma\.V^{-1}\intd{V})=\tr(\sigma\.V^{-1}\intd{F}^T)+\tr(\sigma\.V^{-1}\intd{F})=2\.\tr(\sigma\.V^{-1}\intd{F})\,,
\end{align*}
since $\sigma$ is symmetric and coaxial to $V$. From \eqref{3.7} we therefore obtain
\begin{align}
	\intd{A}=-\det(V)\cdot\tr(\sigma\.V^{-1}\intd{V}).\label{3.8}
\end{align}
Hence, due to \eqref{1.2}, \eqref{1.3} and \eqref{2.4}:
\begin{align*}
	\intd{A}=-\ex{j}\cdot\tr(\sigma\.\intd{L})\,.\tag{3.8*}
\end{align*}
If we substitute this expression into the isothermal relation \eqref{3.5} and use \eqref{eq:repformula}, then it follows:
\begin{align*}
	\ex{j}\cdot [f_1 \tr(\intd{L})+f_2 \tr(L\.\intd{L})+f_3 \tr(L^2\intd{L})]=\pdd{W}{j}\.\intd{j}+\pdd{W}{k}\.\intd{k}+\pdd{W}{l}\.\intd{l}\,.
\end{align*}
Since, by \eqref{2.4},
\begin{align*}
	\intd{j}=\tr(\intd{L})\,,\qquad\intd{k}=2\.\tr(L\.\intd{L})\qquad\text{and}\qquad\intd{l}=3\.\tr(L^2\intd{L})\,,
\end{align*}
we finally conclude that
\begin{align*}
	\ex{j}\.f_1=\pdd{W}{j}\,,\qquad\ex{j}\.f_2=2\.\pdd{W}{k}\,,\qquad\ex{j}\.f_3=3\.\pdd{W}{l}
\end{align*}
and therefore, with \eqref{eq:repformula},
\begin{align}
	\sigma\.\ex{j}=\pdd{W}{j}\cdot\id+2\.\pdd{W}{k}\cdot L+3\.\pdd{W}{l}\cdot L^2\,,\qquad W=W(j,k,l,\Theta)\,.\label{3.9}
\end{align}
Accordingly, from \eqref{3.6} we obtain for the adiabatic law:
\begin{align}
	\sigma\.\ex{j}=\pdd{E}{j}\cdot\id+2\.\pdd{E}{k}\cdot L+3\.\pdd{E}{l}\cdot L^2\,,\qquad E=E(j,k,l,S)\,.\label{3.10}
\end{align}
If we want to omit the introduction of $L$ and use $V$ directly when formulating the law of elasticity, then we appropriately use the following as the invariants of $V$:
\begin{align*}
	I_1=\tr(V)\,,\qquad I_2=\half\.\tr(V^2)\,,\qquad I_3=\det(V)\,.
\end{align*}
Furthermore, according to \eqref{2.7}, \eqref{3.5} and \eqref{3.8}, an analogous computation leads to the law of elasticity in the form
\begin{align}\label{eq:sigmaW}
	\sigma=\pdd{W}{I_3}\cdot\id+\tel{I_3}\cdot\pdd{W}{I_1} \cdot V+\tel{I_3}\cdot\pdd{W}{I_2}\cdot V^2\,,\qquad W=W(I_\nu, \Theta)
\end{align}
and a corresponding formulation with $E(I_1,I_2,I_3,S)$ instead of $W$.
%%%%%%%%%%%%%%%%%%%%%%%%%%%%%%%%%%%%%%%%%%%%%%%%%%%%%%%%%%%%%%%%%%%%%%%%%%%%%%%%%%%
\section{Transition to the deviators}
The introduction of the logarithmic stretch $L$ now proves to be not only appropriate to formulate the law of elasticity as simple as possible, but using $L$ also allows for the decomposition of a deformation into a shape change and volume change by simply taking the deviatoric part, i.e.\ the same approach as for infinitesimal strains, whereas a corresponding decomposition in terms of $V$ is highly inconvenient.
To see this, we decompose the general stretch $V$ into a shape-changing stretch $V_g$ and a volume-changing stretch $V_v$, i.e.\ we demand:
\begin{align}
	V=V_g \cdot V_v=V_v \cdot V_g \quad\text{with}\quad\det V_g=1\quad\text{and}\quad V_v=\beta \cdot\id\quad\text{with}\quad\beta>0\,.\label{4.1}
\end{align}
Obviously, \eqref{4.1} uniquely determines such a decomposition for each $V$ with $\det V>0$; namely, for given $V$,
\begin{align*}
	\beta=\sqrt[3]{\det V}\qquad\text{and}\qquad V_g=\beta^{-1}\cdot V\,.
\end{align*}
Since $V_g$ commutes with $V_v$, we can take the logarithm of \eqref{4.1}:
\begin{align}
	L=L_g+L_v \qquad \text{with} \qquad L_g=\log V_g \qquad \text{and} \qquad L_v=\log V_v\,. \label{4.2}
\end{align}
Then, by \eqref{1.3}, we obtain:
\begin{align*}
	\tr(L_g)=\log(\det V_g)=0\,,\qquad L_v=\log\beta\cdot\id\,,\qquad\tr(L_v)=3\.\log\beta\,.
\end{align*}
If, in general, we denote by $\dev D$ the deviator corresponding to the symmetric matrix $D$, i.e.
\begin{align}
	\dev D=D-\tel{3}\.\tr D\cdot\id\,,\label{4.3}
\end{align}
we can finally write:
\begin{align}
	L_g=\dev L\qquad\text{and}\qquad L_v=\tel{3}\.j\cdot\id\,.\label{4.4}
\end{align}
Thus the change of shape is indeed characterized by the deviator of $L$. For infinitesimal strains we have $L \approx V-\id$, so that $\dev L$ turns into the usual deformation deviator.
\par \noindent
If we now introduce the invariants of $\dev L$:
\begin{alignat}{3}
	y&=\tr((\dev L)^2)\qquad &&\text{and}\qquad &z&=\tr((\dev L)^3)\,,\label{4.5}\\
\intertext{then}
	y&=k-\tel{3}\.j^2\qquad &&\text{and}\qquad &z&=l-j\.k+\frac{2}{9}\.j^3\,.\notag
\end{alignat}
We can use $j$, $y$ and $z$ instead of $j$, $k$ and $l$ as variables. Then $j$ characterizes the change of volume, whereas $y$ and $z$ characterize the change of shape. As one can easily calculate, \eqref{3.9} leads to the formula
\begin{equation}
	\left. \begin{aligned}
	\tel{3}\.\ex{j}\.\tr\sigma&=\pdd{W}{j}\\
	\ex{j}\cdot\dev\sigma&=-y\.\pdd{W}{z}\cdot\id+2\.\pdd{W}{y}\cdot\dev L+3\.\pdd{W}{z}\.(\dev L)^2 \quad
\end{aligned}\right\}
\label{4.6}
\end{equation}
where, in contrast to \eqref{3.9}, $W=W(j,y,z,\Theta)$ now holds.
\par \medskip \noindent
A corresponding formula results from \eqref{3.10}.
\par \medskip \noindent
Without proof, let us remark that $y$ and $z$ cannot take on all possible values independently of each other, but are restricted by the condition
\begin{align*}
	0\leq\frac{z^2}{y^3}\leq\tel{6}\,.
\end{align*}
%%%%%%%%%%%%%%%%%%%%%%%%%%%%%%%%%%%%%%%%%%%%%%%%%%%%%%%%%%%%%%%%%%%%%%%%%%%
\section{Decomposable elasticity laws}
\label{section:richter5}
In the elasticity theory of infinitesimal strains the elastic energy can be interpreted as the sum of the energy of the volume and shape change. Since the change of volume is represented by $j$ and the change of shape is represented by $y$ and $z$, this decomposition is possible for the case of finite strains if and only if
\begin{align}
	\left.\begin{aligned}
		W&=W_{\mathrm{vol}}(j, \Theta)+W_{\mathrm{iso}}(y, z, \Theta)\,,\quad\\
		\text{resp.}\qquad E&=E_{\mathrm{vol}}(j, S)+E_{\mathrm{iso}}(y, z, S)\label{5.1}
	\end{aligned}\right\}
\end{align} 
holds. Then with \eqref{4.6}:
\begin{align*}
	\tel{3}\.\ex{j}\cdot\tr\sigma=\pdd{W_{\mathrm{vol}}}{j}(j, \Theta)\,.
\end{align*}
Thus the average stress depends only on $j$, i.e.\ on the change of volume. If, vice versa, $\tr \sigma$ depends only on $j$, then by \eqref{4.6} we obtain
\begin{align*}
	\frac{\partial^2 W}{\partial j\.\partial y}=\frac{\partial^2 W}{\partial j\.\partial z}=0\,,
\end{align*}
which also leads to the form of $W$ in \eqref{5.1}. Consequently, we can state: \emph{The elastic energy can be decomposed into the energy of change of volume and of change of shape if and only if the mean stress depends only on the change of volume.}
%%%%%%%%%%%%%%%%%%%%%%%%%%%%%%%%%%%%%%%%%%%%%%%%%%%%%%%%%%%%%%
\section{Transition to a new reference temperature}
We referred the deformations to the stress-free state at a certain temperature $\Theta_0$. Now we assume another temperature $\Theta_1$ to be used as initial temperature instead of $\Theta_0$. For $\sigma=0$, the temperature $\Theta_1$ corresponds to a certain deformation $V_1$ with $\log V_1=L_1$. $V_1$ is a scalar multiple of the identity tensor; thus $\dev L_1=0$, $y_1=z_1=0$. Then with \eqref{4.6}:
\begin{align*}
	\pdd{W}{j}(j_1,0,0,\Theta_1)=0\,,
\end{align*}
which leads to the law of thermal expansion:
\begin{align}
	j_1=\varphi(\Theta_1)\,.
\end{align}
Since $\widehat{F}=F\.V_1^{-1}$ is the matrix corresponding to the deformation $F$ with respect to the new initial state, we thus have
$\widehat{V}=V\.V_1^{-1}$, $\widehat{L}=L-L_1$ and hence
\begin{align}
	\widehat{j}=j-j_1\,,\qquad\quad\widehat{y}=y\,,\qquad\quad\widehat{z}=z\,.
\end{align}
In formula \eqref{4.6}, we can now replace $j$ by $\widehat{j}$ if we simultaneously substitute $W$ with
\begin{align}
	\widehat{W}\Big(\widehat{j},y,z,\Theta\Big)=\ex{-j_1}\cdot W\Big(\widehat{j}+j_1,y,z,\Theta\Big)=\ex{-\varphi\big(\Theta_1\big)}\cdot W\Big(\widehat{j}+\varphi\big(\Theta_1\big),y,z,\Theta\Big)\,.
\end{align}
In particular, it follows  
that the decomposition of the elastic energy, which was discussed in Section \ref{section:richter5}, is independent of the choice of the reference temperature.
%%%%%%%%%%%%%%%%%%%%%%%%%%%%%%%%%%%%%%%%%%%%%%%%%%%%%%%%%%%%%%%%%%%%%%%%%%%%%%%
\section{Validity of Hooke's law}
Due to the formulae found previously, one can impose a wide variety of requirements on the law of elasticity, in particular with respect to the dependence on temperature, and verify if these requirements can be satisfied. Let us now consider the question whether the common law by \emph{Hooke} remains valid for finite strains.
\par \noindent
Using the Lam\'{e} constants, \emph{Hooke}'s law takes the form
\begin{align}
	\sigma&=\lambda\cdot\tr(V-\id)\cdot\id+2\.\mu\cdot (V-\id)\label{7.1}\\
\intertext{or}
	\sigma&=(\lambda\cdot I_1-3\.\lambda-2\.\mu)\cdot\id+2\.\mu\cdot V\,.\label{7.2}
\end{align}
It is obvious that \eqref{7.1} is actually derived from the general formula \eqref{3.9} for small $L$.
\par \noindent
In order for the isothermal law of elasticity \eqref{7.2} to remain valid for finite strains, the following equations must be fulfilled according to \eqref{eq:sigmaW}:
\begin{align*}
	\lambda\.I_1-3\.\lambda-2\.\mu=\pdd{W}{I_3}\,,\qquad 2\.\mu\.I_3=\pdd{W}{I_1}\qquad\text{and}\qquad 0=\pdd{W}{I_2}\,.
\end{align*}
This is only possible if $\lambda=2\.\mu$, which corresponds to the Poisson ratio $\nu=\frac{1}{3}$. For all other values of $\nu$, Hooke's law cannot be used for finite strains. Instead, one can use the corresponding logarithmic law
\begin{align}
	\sigma\.\ex{j}=\lambda\.j\cdot\id+2\.\mu\.L\,,\label{7.3}
\end{align}
which, in the isothermal case, corresponds to the decomposable energy
\begin{align*}
	W=\frac{\lambda}{2}\.j^2+\mu\.k=\left(\frac{\lambda}{2}+\frac{\mu}{3}\right)\cdot j^2+\mu\cdot y.
\end{align*}
\let\thefootnote\relax\footnote{Received 2.\ February 1948}
%%%%%%%%%%%%%%%%%%%%%%%%%%%%%%%%%%%%%%%%%%%%%%%%%%%%%%%%%%%%%%%%%%
\newpage
\appendix
\section{List of Symbols}
\footnotesize{\label{appendix:listOfSymbols}
\begin{longtable}{| l | l | l |}
\hline
\textbf{Our notation}	&\textbf{Richter's notation}					&\textbf{meaning}\\
\hline
							&											& \\
$X$, $Y$					&$\textfrak{A}$, $\textfrak{B}$				&arbitrary $3 \times 3$-matrices\\
%$x$, $y$, \dots			&$\textfrak{x}$, $\textfrak{y}$ \dots		&vectors\\												
$X^T$						&$\overline{\textfrak{A}}$					&transpose of $X$\\
$(X)_{ik}$					&$(\textfrak{A})_{ik}$						&entry in the $i$-th row and the $k$-th column of $X$\\
$\det X$					&$|\textfrak{A}|$							&determinant of $X$\\
$\tr X$						&${\textfrak{A}}$							&trace of $X$\\
$\id$						&$\textfrak{E}$								&identity tensor \\
$X^{-1}$					&$\textfrak{A}^{-1}$						&inverse of $X$\\
							&											&\\
$F$							&$\textfrak{A}$								&Jacobian matrix (state of strain)\\
$R$							&$\textfrak{R}$								&pure Euclidean rotation\\
$V$							&$\textfrak{S}$								&pure stretch\\
$\sigma$					&$\textfrak{P}$								&stress tensor (state of stress)\\
$\Theta$					&$\Theta$									&temperature\\
$I_1$, $I_2$, $I_3$			&$I_1$, $I_2$, $I_3$						&invariants of $V$\\
$L$							&$\textfrak{L}$								&logarithmic stretch:\enspace $L=\log V$\\
$j$, $k$, $l$				&$j$, $k$, $l$								&invariants of $L$:\enspace $j=\tr(L)$, $k=\tr(L^2)$, $l=\tr(L^3)$\\
$f_1$, $f_2$, $f_3$			&$f_1$, $f_2$, $f_3$						&coefficient functions\\
$g_1$, $g_2$, $g_3$			&$g_1$, $g_2$, $g_3$						&coefficient functions\\
$X$							&$\textfrak{X}$								&$X=f_1 \cdot \id+f_2 \cdot L+f_3 \cdot L^2$\\
&&\\
$E$							&$u$										&internal energy\\
$S$							&$s$										&entropy\\
$W$							&$f$										&free energy\\
$\intd{A}$					&$\intd{A}$									&differential of the work\\
$e_1$, $e_2$, $e_3$			&$e_1$, $e_2$, $e_3$						&unit vectors in the principal stretch directions of $V$\\
$\sigma_1$, $\sigma_2$, $\sigma_3$	&$\sigma_1$, $\sigma_2$, $\sigma_3$	&components of $\sigma$ in the principal stretch directions of $V$\\
&&\\
$V_g$, $V_v$				&$\textfrak{S}_g$, $\textfrak{S}_v$			&stretch in shape, stretch in volume\\
$\beta$						&$\beta$									&stretch factor of the stretch in volume $V_v$\\
$L_g$, $L_v$				&$\textfrak{L}_g$, $\textfrak{L}_v$			&$L_g=\log V_g$, $L_v=\log V_v$\\
$D$							&$\textfrak{D}$								&arbitrary symmetric matrix\\
$\dev D$					&$\widetilde{\textfrak{D}}$					&common deviator of $D$\\
$y$, $z$					&$y$, $z$									&invariants of $\dev L$:\enspace $y=\tr((\dev L)^2)$, $z=\tr((\dev L)^3)$\\
&&\\
$W_{\mathrm{vol}}$, $E_{\mathrm{vol}}$	&$F$, $U$						&volumetric energies\\
$W_{\mathrm{iso}}$, $E_{\mathrm{iso}}$	&$G$, $V$						&isochoric energies\\
&&\\
$\Theta_0$ resp. $\Theta_1$	&$\Theta_0$ resp. $\Theta_1$				&reference temperatures\\
$\phantom{\Theta}_1$ [index]	&$\phantom{\Theta}_1$				&indicates the correspondence to the temperature $\Theta_1$\\
$\varphi$					&$\varphi$									&logarithmic thermal expansion\\
$\widehat{F}$				&$\textfrak{A}'$							&deformation with respect to the initial state at $\Theta_1$\\
$\widehat{\phantom{F}}$		&$\phantom{\textfrak{A}}'$					&indicates the correspondence to the deformation $\widehat{F}$\\
&&\\
$\lambda$, $\mu$			&$\lambda$, $\mu$							&Lam\'{e} constants\\							
$\nu$						&$m=\tel{\nu}$								&Poisson modulus\\

&&\\
\hline												
\end{longtable}}
\end{document}